\documentclass{article}
\usepackage{amssymb}
\usepackage{amsmath}
\usepackage[all]{xy}

\title{\Large \bf Approximation Presentations of Modules\\
and Homological Conjectures
\thanks{2000 Mathematics Subject
Classification: 16E10, 16E30, 16G10.}
\thanks{Keywords: ${\rm
\mathbb{W}}^{t}$-approximation presentations, property $({\rm
W}^{k})$, homological conjectures, right quasi $k$-Gorenstein
rings.}}
\author{Zhaoyong Huang\thanks{\small \it E-mail: huangzy@nju.edu.cn}\\
{\small \it Department of Mathematics, Nanjing University, Nanjing
210093, People's Republic of China}\\ \\}
\usepackage{amssymb}
\date{}
\begin{document}
\baselineskip=18pt \maketitle

\begin{abstract} In this paper we give a sufficient condition of the
existence of ${\rm \mathbb{W}}^{t}$-approximation presentations. We
also introduce property (W$^{k}$). As an application of the
existence of ${\rm \mathbb{W}}^{t}$-approximation presentations we
give a connection between the finitistic dimension conjecture, the
Auslander-Reiten conjecture and property (W$^{k}$).
\end{abstract}

\vspace{0.5cm}

\centerline{\large \bf 1. Introduction}

\vspace{0.2cm}

In homological algebra and representation theory of algebras, the
following is an important and interesting question, which is
connected with the finitistic dimension conjecture.

{\bf Question:} For a ring $\Lambda$, are the left and right
self-injective dimensions of $\Lambda$ identical?

Zaks [Za] proved that the answer is affirmative for a left and right
noetherian ring if both dimensions are finite. Such rings are called
{\it Gorenstein}.

For a positive integer $k$, Auslander and Reiten in [AR2] initiated
the study of $k$-Gorenstein algebras, which has stimulated several
investigations. They showed that the answer to the question above is
positive in case $\Lambda$ is an artin $\infty$-Gorenstein algebra
(that is, $\Lambda$ is artin $k$-Gorenstein for all $k$). In [AR1],
they also gave the relationship between the question above and the
finitistic dimension conjecture (resp. the contravariant finiteness
of the full subcategory of mod $\Lambda$ consisting of the modules
$M$ with Ext$_{\Lambda}^{i}(M, \Lambda)=0$ for any $i\geq 1$).

It follows from [AR3, Theorem 0.1] and [HN, Theorem 4.1] that the
following conditions are equivalent for a left and right Noetherian
ring $\Lambda$:

(1) For a minimal injective resolution $0 \to \Lambda \to I_0 \to
I_1 \to \cdots \to I_i \to \cdots$ of $\Lambda$ as a right
$\Lambda$-module, the right flat dimension of $I_i$ is at most $i+1$
for any $0 \leq i \leq k-1$.

(2) The strong grade of Ext$_{\Lambda}^{i+1}(M, \Lambda)$ is at
least $i$ for any $M\in$mod $\Lambda$ and $1 \leq i \leq k$.

(3) The grade of Ext$_{\Lambda}^i(N, \Lambda)$ is at least $i$ for
any $N\in$mod $\Lambda ^{\rm op}$ and $1 \leq i \leq k$.

We call a ring {\it right quasi k-Gorenstein} provided it satisfies
one of these equivalent conditions. A ring is called {\it right
quasi} $\infty$-{\it Gorenstein} if it is right quasi $k$-Gorenstein
for all $k$. From [H1] we know that there are right quasi
$k$-Gorenstein rings which are not $k$-Gorenstein; and contrary to
the notion of $k$-Gorenstein, the notion of quasi $k$-Gorenstein is
not left-right symmetric. We showed in [H3] that the answer to the
question above is also positive if $\Lambda$ is an artin right quasi
$\infty$-Gorenstein algebra.

For a ring $\Lambda$ and a positive integer $t$, recall that a
module $M\in$mod $\Lambda$ (resp. mod $\Lambda ^{op}$) is called a
${\rm W}^{t}$-{\it module} if Ext$_{\Lambda}^{i}(M, \Lambda)=0$ for
any $1\leq i \leq t$. We remark that $i=0$ is not required in this
definition. $M$ is called a ${\rm W}^{\infty}$-{\it module} if it is
a ${\rm W}^{t}$-module for all $t$. Jans called in [J2] a module
$M\in$mod $\Lambda$ (resp. mod $\Lambda ^{op}$) a W-{\it module} if
Ext$_{\Lambda}^{1}(M, \Lambda)=0$. This is the motivation for us to
give the above definition of ${\rm W}^{t}$-modules. We use
$\mathbb{W}^{t}$  (resp. $\mathbb{W}^{\infty}$) to denote the full
subcategory of mod $\Lambda$ consisting of ${\rm W}^{t}$-modules
(resp. ${\rm W}^{\infty}$-modules). We call an exact sequence $0 \to
K_{t}(M) \to E_{t}(M) \to M \to 0$ a $\mathbb{W}^{t}$-{\it
approximation presentation} of $M$ if $E_{t}(M)\to M$ is a right
$\mathbb{W}^{t}$-approximation of $M$ and the projective dimension
of $K_{t}(M)$ is at most $t-1$ (see [H1]).

One of the main results in [H1] is that if $\Lambda$ is a right
quasi $k$-Gorenstein algebra then every module in mod $\Lambda$ has
a $\mathbb{W}^{t}$-approximation presentation for any $1\leq t \leq
k$. In Section 3 we give a sufficient condition for the existence of
${\rm \mathbb{W}}^{t}$-approximation presentations. Let $\Lambda$ be
a left and right noetherian ring and $k$ a positive integer. For any
$1\leq t\leq k$ we show that a module $M$ in mod $\Lambda$ has a
${\rm \mathbb{W}}^{t}$-approximation presentation if the strong
grade of Ext$_{\Lambda}^{i+1}(M, \Lambda)$ is at least $i$ for any
$1\leq i \leq k-1$. This improves the main result in [H1] and [AB,
Proposition 2.21]. We then study the homological finiteness of the
full subcategory of mod $\Lambda$ (resp. mod $\Lambda ^{op}$)
consisting of the modules with projective dimension at most $k$. In
particular, we show that, over an artin right quasi $k$-Gorenstein
algebra $\Lambda$, such subcategory is functorially finite in mod
$\Lambda ^{op}$.

In [AR2] Auslander and Reiten posed the conjecture ({\bf ARC}): An
artin algebra is Gorenstein if it is $\infty$-Gorenstein. The famous
Nakayama conjecture is a special case of this conjecture. We
introduce in Section 2 {\it property} (W$^{k}$): Each W$^{k}$-module
in mod $\Lambda$ is torsionless. We then give some equivalent
conditions for property (W$^{k}$), which we subsequently use to
prove the following result: If the strong Nakayama conjecture holds
true for $\Lambda$, and the left self-injective dimension of
$\Lambda$ is at most one, then the right self-injective dimension of
$\Lambda$ is at most one as well.

As an application of the results obtained in Section 3 we show in
Section 4 the validity of {\bf ARC} is equivalent to property
(W$^{k}$), and we give the relationship with the finitistic
dimension. Let $\Lambda$ be an artin right quasi $\infty$-Gorenstein
algebra and $k$ a non-negative integer. We show that $\Lambda$ is
Gorenstein with self-injective dimension at most $k$ if and only if
$\Lambda$ has property (W$^{k})$, and that the difference between
the self-injective dimension and the finitistic dimension of
$\Lambda$ is at most one. We then conclude that in order to verify
Nakayama conjecture (for any artin algebra) it suffices to verify
the finitistic dimension conjecture for (right quasi)
$\infty$-Gorenstein algebras.

Let $M\in$ mod $\Lambda$ and $n$ a non-negative integer. If $M$
admits a resolution (of finite length) $0 \to X_{n} \to \cdots \to
X_{1} \to X_{0} \to M \to 0$ with all $X_{i}\in
{\mathbb{W}^{\infty}}$, then set
$\mathbb{W}^{\infty}$-dim$_{\Lambda}M=$inf$\{ n|$ there is an exact
sequence $0 \to X_{n} \to \cdots \to X_{1} \to X_{0} \to M \to 0$
with all $X_{i}\in {\mathbb{W}^{\infty}} \}$. If no such a
resolution exists, set
$\mathbb{W}^{\infty}$-dim$_{\Lambda}M=\infty$. We call
$\mathbb{W}^{\infty}$-dim$_{\Lambda}M$ the {\it left orthogonal
dimension} of $M$ (see [H2]).

In Section 5 we show that if $\Lambda$ has property (W$^{\infty}$)
then for each $M\in$mod $\Lambda$ the left orthogonal dimension of
$M$ and its Gorenstein dimension are identical, which yields that
$\Lambda$ has property (W$^{\infty}$) if and only if each
W$^{\infty}$-module in mod $\Lambda$ is reflexive. We also show that
{\bf ARC} is true if and only if $\mathbb{W}^{\infty}$ is
contravariantly finite and the global dimension of $\Lambda$
relative to $\mathbb{W}^{\infty}$ is finite.

According to the results obtained in the former sections, we pose in
the final section two conjectures: (1) Any artin algebra has
property (W$^{\infty}$). (2) An artin algebra is Gorenstein if it is
right quasi $\infty$-Gorenstein. The latter conjecture is clearly a
generalized version of {\bf ARC}.

{\bf Definitions and notations}

In the following we give some definitions and notations which are
often used in this paper.

For a ring $\Lambda$, we use mod $\Lambda$ to denote the category of
finitely generated left $\Lambda$-modules and
$\mathcal{P}^{k}(\Lambda)$ (resp. $\mathcal{P}^{\infty}(\Lambda)$)
to denote the full subcategory of mod $\Lambda$ consisting of the
modules with projective dimension at most $k$ (resp. the modules
with finite projective dimension). For a left (resp. right)
$\Lambda$-module $M$, l.pd$_{\Lambda}M$ and l.id$_{\Lambda}M$ (resp.
r.pd$_{\Lambda}M$ and r.id$_{\Lambda}M$) are denoted the left
projective dimension and the left injective dimension (resp. the
right projective dimension and the right injective dimension) of
$M$, respectively. In addition, we recall from [ASm1] the following
definition.

Assume that $\mathcal{C}\supset \mathcal{D}$ are full subcategories
of mod $\Lambda$ (resp. mod $\Lambda ^{op}$) and $C\in \mathcal{C}$,
$D\in$add$\mathcal{D}$, where add$\mathcal{D}$ is the full
subcategory of mod $\Lambda$ (resp. mod $\Lambda ^{op}$) consisting
of all $\Lambda$-modules (resp. $\Lambda ^{op}$-modules) isomorphic
to summands of finite direct sums of modules in $\mathcal{D}$. A
morphism $D\to C$ is said to be a right $\mathcal{D}$-{\it
approximation} of $C$ if Hom$_{\Lambda}(X, D)\to$Hom$_{\Lambda}(X,
C)\to 0$ is exact for all $X\in$add$\mathcal{D}$. The subcategory
$\mathcal{D}$ is said to be {\it contravariantly finite} in
$\mathcal{C}$ if every module $C$ in $\mathcal{C}$ has a right
$\mathcal{D}$-approximation. Dually, we define the notions of left
$\mathcal{D}$-{\it approximation} and {\it covariantly finite}. The
subcategory $\mathcal{D}$ is said to be {\it functorially finite} in
$\mathcal{C}$ if it is both contravariantly finite and covariantly
finite in $\mathcal{C}$. The notions of contravariantly finite
subcategories, covariantly finite subcategories and functorially
finite subcategories are referred to as {\it homologically finite
subcategories}.

We also list some famous homological conjectures.

{\bf Finitistic Dimension Conjecture (FDC):} fin.dim$\Lambda
$=sup$\{$l.pd$_{\Lambda}M|M\in$ mod $\Lambda$ and l.pd$_{\Lambda}M$
is finite$\}$ is finite for any artin algebra $\Lambda$.

A brief history and some recent development of {\bf FDC} were given
in [Zi].

{\bf Auslander-Reiten Conjecture (ARC):} Every artin
$\infty$-Gorenstein algebra is Gorenstein.

{\bf Nakayama Conjecture (NC):} An artin algebra $\Lambda$ is
self-injective if each term in a minimal injective resolution of
$\Lambda$ as a right $\Lambda$-module is projective.

Let us recall the relationship between the conjectures mentioned
above. It is shown in [Y] that {\bf FDC} implies {\bf NC} and in
[AR2] that {\bf ARC} implies {\bf NC}.

Throughout this paper, unless stated otherwise, $\Lambda$ is a left
and right noetherian ring.

\vspace{0.5cm}

\centerline{\large \bf 2. Property (W$^{k}$)}

\vspace{0.2cm}

In this section we introduce property (W$^{k}$) and then give some
characterizations and applications.

Let $A\in$mod $\Lambda$ (resp. mod $\Lambda ^{op}$) and $\sigma
_{A}: A \to A^{**}$ defined via $\sigma _{A}(x)(f)=f(x)$ for any
$x\in A$ and $f\in A^{*}$. $A$ is called {\it torsionless} if
$\sigma _{A}$ is a monomorphism; and $A$ is called {\it reflexive}
if $\sigma _{A}$ is an isomorphism. Let $P_1 \to P_0 \to A \to 0$ be
a projective resolution of $A$ in mod $\Lambda$ (resp. mod $\Lambda
^{op}$). For a positive integer $k$, $A$ is called $k$-{\it
torsionfree} if Ext$_{\Lambda}^i({\rm Tr}A, \Lambda)=0$ for any $1
\leq i \leq k$, where ${\rm Tr}A={\rm Coker}(P_0^* \to P_1^*)$ is
the {\it transpose} of $A$ (see [AB] or [AR3]). We remark that it is
known that a module in mod $\Lambda$ is torsionless (resp.
reflexive) if and only if it is 1-torsionfree (resp. 2-torsionfree)
(see [AB]).

The following lemma, which is of independent interest, is useful for
the rest of this paper.

\vspace{0.2cm}

{\bf Lemma 2.1} {\it For a positive integer} $k$, {\it the
following statements are equivalent}.

(1) {\it Each} ${\rm W}^{k}$-{\it module in} mod $\Lambda$ {\it is
torsionless}.

(2) {\it Each} ${\rm W}^{k}$-{\it module in} mod $\Lambda$ {\it is
reflexive}.

(3) {\it Each} ${\rm W}^{k}$-{\it module in} mod $\Lambda ^{op}$
{\it is} ${\rm W}^{\infty}$-{\it module}.

\vspace{0.2cm}

{\it Proof.} $(2)\Rightarrow (1)$ It is trivial.

$(1)\Rightarrow (3)$ The proof is essentially the same as that of
$(1)\Rightarrow (2)$ in [HT, Lemma 3.3]. Let $Y$ be a ${\rm
W}^{k}$-module in mod $\Lambda ^{op}$ and $P_{k+1} \buildrel
{d_{k+1}} \over \longrightarrow P_{k} \buildrel {d_{k}} \over \to
\cdots \buildrel {d_{1}} \over \to P_{0} \to Y \to 0$ a projective
resolution of $Y$ in mod $\Lambda ^{op}$. It is not difficult to
verify that Coker$d_{k+1}^{*}$ is a ${\rm W}^{k}$-module in mod
$\Lambda$. By (1) Coker$d_{k+1}^{*}$ is torsionless. Then by [HT,
Lemma 2.1], Ext$_{\Lambda}^{1}({\rm Coker}d_{k+1}, \Lambda)=0$ and
Ext$_{\Lambda}^{k+1}(Y, \Lambda)(\cong$Ext$_{\Lambda}^{1}({\rm
Coker}d_{k+1}, \Lambda))=0$.

Since $0\to {\rm Im}d_{1} \to P_{0} \to Y \to 0$ is exact,
Ext$_{\Lambda}^{i}({\rm Im}d_{1}, \Lambda)=0$ for any $1 \leq i
\leq k$. Repeating the above argument we have
Ext$_{\Lambda}^{k+1}({\rm Im}d_{1}, \Lambda)=0$ and thus
Ext$_{\Lambda}^{k+2}(Y, \Lambda)=0$. Continuing this procedure, we
get our conclusion.

$(3)\Rightarrow (2)$ By [HT, Lemma 3.3]. \hfill{$\blacksquare$}

\vspace{0.2cm}

{\bf Definition 2.2} We say that $\Lambda$ has {\it property}
(W$^{k}$) if the condition (1) of Lemma 2.1 is satisfied for
$\Lambda$. If $k$ is infinite, we then say that $\Lambda$ has {\it
property} (W$^{\infty}$). We say that $\Lambda$ has {\it property}
(W$^{k})^{op}$ (resp. (W$^{\infty})^{op}$) if $\Lambda ^{op}$ has
property (W$^{k})$ (resp. (W$^{\infty}$)).

\vspace{0.2cm}

{\it Remark.} (1) It is trivial that if $\Lambda$ has property
(W$^{k}$) (resp. (W$^{k})^{op}$) then $\Lambda$ has property
(W$^{n}$) (resp. (W$^{n})^{op}$) for any $n>k$ and property
(W$^{\infty}$) (resp. (W$^{\infty})^{op}$).

(2) In Section 5 we will show that (1) and (2) in Lemma 2.1 are also
equivalent when $k$ is infinite. That is, $\Lambda$ has property
(W$^{\infty})$ if and only if each W$^{\infty}$-module in mod
$\Lambda$ is reflexive (see Corollary 5.2).

\vspace{0.2cm}

Let $N \in$mod $\Lambda^{op}$ and
$$0 \to N \buildrel {\delta _0} \over \longrightarrow E_{0}
\buildrel {\delta _1} \over \longrightarrow E_{1} \buildrel {\delta
_2} \over \longrightarrow \cdots \buildrel {\delta _i} \over
\longrightarrow E_{i} \buildrel {\delta _{i+1}} \over
\longrightarrow \cdots$$ be an injective resolution of $N$. Recall
from [CF] that an injective resolution as above is called {\it
ultimately closed} at $n$ (where $n$ is a positive integer) if
Im$\delta _n=\bigoplus _{j=0}^mW_j$, where each $W_j$ is a direct
summand of Im$\delta _{i_j}$ with $i_j <n$.

\vspace{0.2cm}

{\bf Corollary 2.3} $\Lambda$ {\it has property} (W$^{k}$) {\it if}

(1) r.id$_{\Lambda}\Lambda \leq k$, {\it or}

(2) $\Lambda _{\Lambda}$ {\it has a ultimately closed injective
resolution at} $k$.

\vspace{0.2cm}

{\it Proof.} Our conclusions follow from [HT, Theorems 2.2 and 2.4],
respectively. \hfill{$\blacksquare$}

\vspace{0.2cm}

{\bf Lemma 2.4} ([AB, Proposition 2.6]) {\it Let} $A \in$mod
$\Lambda$ ({\it resp.} mod $\Lambda^{op}$). {\it Then we have the
following exact sequences:}
$$0 \to{\rm Ext}^{1}_{\Lambda}({\rm Tr}A, \Lambda) \to A \buildrel {\sigma
_{A}} \over \longrightarrow A^{**} \to {\rm Ext}^{2}_{\Lambda}({\rm
Tr}A, \Lambda)\to 0 \eqno{\textcircled{\small{1}}}$$
$$0 \to{\rm Ext}^{1}_{\Lambda}(A, \Lambda) \to {\rm Tr}A \buildrel {\sigma
_{{\rm Tr}A}} \over \longrightarrow ({\rm Tr}A)^{**} \to {\rm
Ext}^{2}_{\Lambda}(A, \Lambda)\to 0 \eqno{\textcircled{\small{2}}}$$

\vspace{0.2cm}

We observe that in case $k=1$ the converse of Corollary 2.3(1) holds
true. That is, we have

\vspace{0.2cm}

{\bf Corollary 2.5} ([Ba, Theorem 3.3] and [AR1, Proposition 2.2])
{\it The following statements are equivalent.}

(1) r.id$_{\Lambda}\Lambda \leq 1$.

(2) $\Lambda$ {\it has property} (W$^{1}$).

(3) {\it Each torsionless module in} mod $\Lambda$ {\it is
reflexive.}

(4) {\it Each torsionless module in} mod $\Lambda ^{op}$ {\it is
a} ${\rm W}^{1}$-{\it module.}

\vspace{0.2cm}

{\it Proof.} The equivalence of (1), (3) and (4) follows from [Ba,
Theorem 3.3]. By Corollary 2.3 we have (1) implies (2).

$(2)\Rightarrow (1)$ Let $A$ be torsionless in mod $\Lambda ^{op}$.
Then there are exact sequences \textcircled{\small{1}} and
\textcircled{\small{2}} as in Lemma 2.4.
 Since $A$ is
torsionless, by the the exactness of the sequence
\textcircled{\small{1}} we have ${\rm Ext}^{1}_{\Lambda}({\rm Tr}A,
\Lambda)=0$, and then by the assumption (2) ${\rm Tr}A$ is
torsionless. So from the exact sequence \textcircled{\small{2}}, we
have ${\rm Ext}^{1}_{\Lambda}(A, \Lambda)=0$. Then it is easy to see
that r.id$_{\Lambda}\Lambda \leq 1$. \hfill{$\blacksquare$}

\vspace{0.2cm}

Following Colby and Fuller [CF], we say the strong Nakayama
conjecture ({\bf SNC}) is true for $\Lambda$ if for any $M\in$ mod
$\Lambda$ the condition Ext$_{\Lambda}^{i}(M, \Lambda)=0$ for all
$i\geq 0$ implies $M=0$. By the proof of [Y, Theorem 3.4.3], for any
artin algebra we have {\bf FDC}$\Rightarrow${\bf SNC}.

\vspace{0.2cm}

{\bf Lemma 2.6} {\it Suppose that} {\bf SNC} {\it holds true for}
$\Lambda$. {\it If each torsionless module in} mod $\Lambda ^{op}$
{\it is reflexive, then each torsionless module in} mod $\Lambda
^{op}$ {\it is a} ${\rm W}^{1}$-{\it module.}

\vspace{0.2cm}

{\it Proof.} Assume that $A$ is a torsionless module in mod $\Lambda
^{op}$. Then $A$ is reflexive by assumption. By [J1, Theorem 1.1],
there are a torsionless module $B$ in mod $\Lambda$ and an exact
sequence:
$$0 \to B \buildrel {\sigma _{B}} \over \longrightarrow B^{**} \to
{\rm Ext}^{1}_{\Lambda}(A, \Lambda) \to 0
\eqno{\textcircled{\small{3}}}$$ Since $B^{*}$ is torsionless in mod
$\Lambda ^{op}$, $B^{*}$ is reflexive (by assumption) and $\sigma
_{B^{*}}$ is an isomorphism. On the other hand, by [AF, Proposition
20.14] we have $\sigma _{B}^{*}\sigma _{B^{*}}=1_{B^{*}}$, it
follows that $\sigma _{B}^{*}$ is also an isomorphism.

By assumption and the dual version of Corollary 2.5,
l.id$_{\Lambda}\Lambda \leq 1$. Because $B^{**} \in$mod $\Lambda$ is
torsionless, there is a module $C\in$mod $\Lambda$ such that
$B^{**}$ is a 1-syzygy of $C$ and ${\rm Ext}^{1}_{\Lambda}(B^{**},
\Lambda)\cong {\rm Ext}^{2}_{\Lambda}(C, \Lambda)=0$.

From the exact sequence \textcircled{\small{3}} we get a long exact
sequence:
$$0\to [{\rm Ext}^{1}_{\Lambda}(A, \Lambda)]^{*} \to B^{***}
\buildrel {\sigma _{B}^{*}} \over \longrightarrow B^{*} \to {\rm
Ext}^{1}_{\Lambda}({\rm Ext}^{1}_{\Lambda}(A, \Lambda), \Lambda) \to
{\rm Ext}^{1}_{\Lambda}(B^{**}, \Lambda)=0$$ So $[{\rm
Ext}^{1}_{\Lambda}(A, \Lambda)]^{*}=0$ and ${\rm
Ext}^{1}_{\Lambda}(A, \Lambda)$ is a ${\rm W}^{1}$-module in mod
$\Lambda$. Because l.id$_{\Lambda}\Lambda \leq 1$,
Ext$_{\Lambda}^{i}({\rm Ext}^{1}_{\Lambda}(A, \Lambda),
\Lambda)$\linebreak $=0$ for any $i\geq 0$. Since {\bf SNC} holds
true for $\Lambda$, ${\rm Ext}^{1}_{\Lambda}(A, \Lambda)=0$.
\hfill{\hfill{$\blacksquare$}}

\vspace{0.2cm}

We now give a partial answer to the question mentioned in the
Introduction as follows.

\vspace{0.2cm}

{\bf Theorem 2.7} {\it If} {\bf SNC} {\it holds true for} $\Lambda$,
{\it then} l.id$_{\Lambda}\Lambda \leq 1$ {\it implies}
r.id$_{\Lambda}\Lambda \leq 1$. {\it In particular, if} {\bf SNC}
{\it is always true, then} l.id$_{\Lambda}\Lambda \leq 1$ {\it if
and only if} r.id$_{\Lambda}\Lambda \leq 1$.

\vspace{0.2cm}

{\it Proof.} Assume that {\bf SNC} holds true for $\Lambda$ and
l.id$_{\Lambda}\Lambda \leq 1$. Then by the dual version of
Corollary 2.5, each torsionless module in mod $\Lambda ^{op}$ is
reflexive. It follows from Lemma 2.6 and Corollary 2.5 that each
torsionless module in mod $\Lambda ^{op}$ is a ${\rm W}^{1}$-module
and r.id$_{\Lambda}\Lambda \leq 1$. \hfill{$\blacksquare$}

\vspace{0.2cm}

It would be interesting to know whether the general case of Theorem
2.7 holds true. That is, if {\bf SNC} holds true for $\Lambda$, then
for any positive integer $k$, does l.id$_{\Lambda}\Lambda \leq k$
imply r.id$_{\Lambda}\Lambda \leq k$? If the answer is affirmative,
then we get that {\bf SNC}$\Rightarrow$the Gorenstein Symmetric
Conjecture ({\bf GSC}). The latter conjecture states that the left
and right self-injective dimensions are identical for any artin
algebra.

For a positive integer $t$, we use $(\mathbb{W}^{t})^{op}$ (resp.
$(\mathbb{W}^{\infty})^{op}$) to denote the subcategory of mod
$\Lambda ^{op}$ consisting of W$^{t}$ (resp. W$^{\infty}$)-modules.

\vspace{0.2cm}

{\bf Proposition 2.8}  (1) {\bf SNC} {\it is true for} $\Lambda$
{\it if} $\Lambda$ {\it has property} (W$^{\infty}$).

(2) {\it If} l.id$_{\Lambda}\Lambda \leq k$ {\it (where} $k\leq
2${\it ), then} $\Lambda$ {\it has property} (W$^{\infty}$) {\it if
and only if} r.id$_{\Lambda}\Lambda \leq k$.

\vspace{0.2cm}

{\it Proof.} (1) directly follows from the definition of property
($W^{\infty}$).

(2) By Corollary 2.3, we get the sufficiency. In the following, we
prove the necessity.

The case for $k=0$ is trivial, and the case for $k=1$ follows from
(1) and Theorem 2.7.

Now suppose that $k=2$ and $A\in$mod $\Lambda ^{op}$ and $P_1
\buildrel {f} \over \to P_0 \to A \to 0$ is a projective resolution
of $A$ in mod $\Lambda ^{op}$. Then we get an exact sequence in mod
$\Lambda$:
$$0 \to A^* \to P_0^* \buildrel {f^*} \over
\to P_1^* \to {\rm Tr}A \to 0.$$ It is not difficult to see that
Ker$f\cong ({\rm Tr}A)^*$. Since l.id$_{\Lambda}\Lambda \leq 2$,
${\rm Ext}^i_{\Lambda}(A^*, \Lambda) \cong {\rm
Ext}^{i+2}_{\Lambda}({\rm Tr}A, \Lambda)=0$ for any $i\geq 1$ and
$A^*\in \mathbb{W}^{\infty}$. By assumption, $\Lambda$ has property
(W$^{\infty})$, so $A^*$ is reflexive. Because $({\rm Tr}A)^*\in$mod
$\Lambda ^{op}$, we have that $({\rm Tr}A)^{**}$ is reflexive by the
above argument. Then $({\rm Tr}A)^{***}$ is also reflexive. On the
other hand, by [AF, Proposition 20.14] we have $\sigma _{{\rm
Tr}A}^{*}\sigma _{({\rm Tr}A)^{*}}=1_{({\rm Tr A})^{*}}$, it follows
that $\sigma _{{\rm Tr}A}^{*}$ is a split epimorphism. Then, by
applying the functor Hom$_{\Lambda}(\ , \Lambda)$ to the exact
sequence \textcircled{\small{2}} in Lemma 2.4, we have that $({\rm
Tr A})^{*}$ is a direct summand of $({\rm Tr A})^{***}$ and hence
$({\rm Tr A})^{*}$ is reflexive. So, by assumption and Lemma 2.4, we
have that ${\rm Tr}({\rm Tr A})^{*} \in
\mathbb{W}^{2}(=\mathbb{W}^{\infty})$ and ${\rm Tr}({\rm Tr A})^{*}$
is reflexive. Then $({\rm Tr A})^{*}\in (\mathbb{W}^{2})^{op}$ again
by Lemma 2.4. Since ${\rm Ext}^{i+2}_{\Lambda}(A, \Lambda) \cong
{\rm Ext}^i_{\Lambda}({\rm Ker}f, \Lambda) \cong {\rm
Ext}^i_{\Lambda}(({\rm Tr}A)^*, \Lambda)$ for any $i\geq 1$, ${\rm
Ext}^3_{\Lambda}(A, \Lambda)=0={\rm Ext}^4_{\Lambda}(A, \Lambda)$.
Thus we conclude that r.id$_{\Lambda}\Lambda \leq 2$.
\hfill{$\blacksquare$}

\vspace{0.2cm}

We wonder whether the assertion in Proposition 2.8(2) holds true
when $k$ is any positive integer. The answer is affirmative when
$\Lambda$ is an artin algebra (see [AHT, Theorem 3.10]).

\vspace{0.5cm}

\centerline{\large \bf 3. The existence of ${\rm
\mathbb{W}}^{t}$-approximation presentations}

\vspace{0.2cm}

Let $M$ be in mod $\Lambda$ (resp. mod $\Lambda ^{op}$) and $i$ a
non-negative integer. We say that the {\it grade} of $M$, written as
grade$M$, is at least $i$ if Ext$_{\Lambda}^j(M, \Lambda)=0$ for any
$0 \leq j < i$. We say that the {\it strong grade} of $M$, written
as s.grade$M$, is at least $i$ if grade$X \geq i$ for each submodule
$X$ of $M$ (see [AR3]). We showed in [H1] that every module in mod
$\Lambda$ has a $\mathbb{W}^{t}$-approximation presentation for any
$1\leq t \leq k$ if $\Lambda$ is a right quasi $k$-Gorenstein
algebra. Actually, the argument we use in proving [H1, Theorem 1]
proves the following more general result (or c.f. [AB, Proposition
2.21]).

\vspace{0.2cm}

{\bf Theorem 3.1} {\it Let} $M$ {\it be in} mod $\Lambda$ ({\it
resp}. mod $\Lambda ^{op}$). {\it If} gradeExt$_{\Lambda}^{t}(M,
\Lambda)\geq t$ {\it for any} $1\leq t \leq k$, {\it then} $M$ {\it
has a} $\mathbb{W}^{t}$-{\it approximation presentation for any}
$1\leq t \leq k$.

\vspace{0.2cm}

In this section we develop this result and show that a module
$M\in$mod $\Lambda$ has a $\mathbb{W}^{t}$-approximation
presentation for any $1\leq t \leq k$ if
s.gradeExt$^{t+1}_{\Lambda}(M, \Lambda)\geq t$ for any $1\leq t \leq
k-1$.

\vspace{0.2cm}

Let $M$ be in mod $\Lambda$ and $P \buildrel {f} \over \to M^{*} \to
0$ an exact sequence in mod $\Lambda ^{op}$ with $P$ projective, and
let $h$ be the composition: $M \buildrel {\sigma _{M}} \over
\longrightarrow M^{**} \buildrel {f^{*}} \over \hookrightarrow
P^{*}$. Set $\mathcal{P}^{0}(\Lambda)=\{$projective modules in mod
$\Lambda \}$. From the proof of [H2, Lemma 1] we have the following

\vspace{0.2cm}

{\bf Lemma 3.2}  $h$ {\it is a left} $\mathcal{P}^{0}(\Lambda)$-{\it
approximation of} $M$.

\vspace{0.1cm}

{\bf Theorem 3.3} {\it Let} $M$ {\it be in} mod $\Lambda$ {\it and}
$k$ {\it a positive integer.  If} s.gradeExt$^{t+1}_{\Lambda}(M,
\Lambda)\geq t$ {\it for any} $1\leq t \leq k-1$, {\it then} $M$
{\it has a} ${\rm \mathbb{W}}^{t}$-{\it approximation presentation
for any} $1\leq t \leq k$.

\vspace{0.1cm}

{\it Proof.} We proceed by induction on $k$. The case $k=1$ follows
from [T, Lemma 6.9]. Now suppose that $k \geq 2$ and a module $M
\in$ mod $\Lambda$ satisfies s.gradeExt$^{t+1}_{\Lambda}(M,
\Lambda)\geq t$ for any $1\leq t \leq k-1$.

Let $0 \to L \buildrel {f} \over \to P \buildrel {g} \over \to M \to
0$ be an exact sequence in mod $\Lambda$ with $P$ projective. Then
Ext$^{t}_{\Lambda}(L, \Lambda)\cong$ Ext$^{t+1}_{\Lambda}(M,
\Lambda)$ for any $t\geq 1$ and s.gradeExt$^{t}_{\Lambda}(L,
\Lambda)\geq t$ for any $1\leq t \leq k-1$. By inductive hypothesis,
$L$ has a ${\rm \mathbb{W}}^{k-1}$-approximation presentation:
$$0 \to K_{k-1}(L) \buildrel {\alpha} \over \to E_{k-1}(L)
\buildrel {\beta} \over \to L \to 0,$$ where
l.pd$_{\Lambda}K_{k-1}(L)\leq k-2$ and $E_{k-1}(L)\in
\mathbb{W}^{k-1}$. So Ext$^{t}_{\Lambda}(K_{k-1}(L), \Lambda)
\cong$Ext$^{t+1}_{\Lambda}(L, \Lambda)$ for any $1\leq t \leq k-2$
and hence s.gradeExt$^{t}_{\Lambda}(K_{k-1}(L), \Lambda)\geq t+1$
for any $1\leq t \leq k-2$.

Notice that $K_{k-1}(L)$ is torsionless by [AB, Proposition 3.17]
(we remark that $K_{k-1}(L)$ is trivially torsionless even for
$k=1$). On the other hand, $L$ is torsionless since it is a
submodule of the projective module $P$. Then by [AR3, Theorem 1.1],
$E_{k-1}(L)$ is torsionless.

Let $Q \buildrel {h_{1}} \over \to [E_{k-1}(L)]^{*} \to 0$ be an
exact sequence in mod $\Lambda ^{op}$ with $Q$ projective. Then $0
\to [E_{k-1}(L)]^{**}\buildrel {h_{1}^{*}} \over \to Q^{*}$ is exact
in mod $\Lambda$. Let $h$ be the composition: $E_{k-1}(L) \buildrel
{\sigma _{E_{k-1}(L)}} \over \longrightarrow
[E_{k-1}(L)]^{**}\buildrel {h_{1}^{*}} \over \to Q^{*}$. Since
$E_{k-1}(L)$ is torsionless, $\sigma _{E_{k-1}(L)}$ and $h$ are
monomorphisms. By Lemma 3.2, $h$ is a left
$\mathcal{P}^{0}(\Lambda)$-approximation of $E_{k-1}(L)$. So there
is a homomorphism $\delta :Q^{*} \to P$ such that $\delta h=f\beta$
and hence we have the following commutative diagram with exact rows
and columns:

$$\xymatrix{& 0 \ar@{-->}[d] & 0 \ar[d] & 0 \ar[d]& &\\
0 \ar[r] & K_{k-1}(L) \ar@{-->}[d]^{\gamma} \ar[r]^{\alpha} &
E_{k-1}(L) \ar[d]^{\binom h 0}
\ar[r]^{\beta} & L \ar[d]^{f} \ar[r] & 0\\
0 \ar[r] & Q^{*} \ar[r] & Q^{*}\oplus P
\ar[r]^{(\delta , \ 0)} & P \ar[d]^{g} \ar[r] & 0\\
& & & M \ar[d] & &\\
& & & 0 & &}
$$
where $\gamma$ is an induced homomorphism. Put
$K_{k}(M)=$Coker$\gamma$ and $E_{k}(M)=$Coker${\binom h 0}$. By
the Snake Lemma we have an exact sequence:
$$0 \to K_{k}(M) \to E_{k}(M) \to M \to 0
\eqno{\textcircled{\small{4}}}$$ By the exactness of $0 \to
K_{k-1}(L) \buildrel {\gamma} \over \to Q^{*} \to K_{k}(M) \to 0$ we
have l.pd$_{\Lambda}K_{k}(M)\leq k-1$. On the other hand, by the
exactness of $0 \to E_{k-1}(L) \buildrel {\binom h 0} \over \to
Q^{*}\oplus P \to E_{k}(M) \to 0$ we have
Ext$_{\Lambda}^{t}(E_{k-1}(L), \Lambda)\cong$Ext$_{\Lambda}^{t+1}$
$(E_{k}(M), \Lambda)$ for any $t\geq 1$, which implies that
Ext$_{\Lambda}^{t}(E_{k}(M), \Lambda)=0$ for any $2 \leq t \leq k$.
In addition, It is easy to see that ${\binom h 0}$ is also a left
$\mathcal{P}^{0}(\Lambda)$-approximation of $E_{k-1}(L)$, thus
Ext$_{\Lambda}^{1}(E_{k}(M), \Lambda)=0$ and we conclude that
Ext$_{\Lambda}^{t}(E_{k}(M), \Lambda)=0$ for any $1 \leq t \leq k$.
So $E_{k}(M)\in {\rm \mathbb{W}}^{t}$ and the exact sequence
\textcircled{\small{4}} is as required. \hfill{$\blacksquare$}

\vspace{0.2cm}

Let
$$0\to \Lambda _{\Lambda}\to I_{0} \to I_{1} \to \cdots \to I_{i}
\to \cdots$$ be a minimal injective resolution of $\Lambda$ as a
right $\Lambda$-module.

Recall that $\Lambda$ is called $k$-{\it Gorenstein} if the right
flat dimension of $I_{i}$ is at most $i$ for any $0\leq i \leq k-1$,
and $\Lambda$ is called $\infty$-{\it Gorenstein} if $\Lambda$ is
$k$-Gorenstein for all $k$. It is well known that the notion of
$k$-Gorenstein rings (resp. $\infty$-Gorenstein rings) is left-right
symmetric (see [FGR, Auslander's Theorem 3.7]).

\vspace{0.2cm}

{\bf Definition 3.4} $\Lambda$ is called {\it right quasi} $k$-{\it
Gorenstein} if the right flat dimension of $I_{i}$ is at most $i+1$
for any $0\leq i \leq k-1$, and $\Lambda$ is called {\it right
quasi} $\infty$-{\it Gorenstein} if it is right quasi $k$-Gorenstein
for all $k$.

\vspace{0.2cm}

By [AR3, Theorem 0.1] and [HN, Theorem 4.1], we have that $\Lambda$
is right quasi $k$-Gorenstein if and only if
s.gradeExt$_{\Lambda}^{i+1}(M, \Lambda)\geq i$ for any $M\in$mod
$\Lambda$ and $1 \leq i \leq k$, if and only if
gradeExt$_{\Lambda}^i(N, \Lambda)\geq i$ for any $N\in$mod $\Lambda
^{\rm op}$ and $1 \leq i \leq k$.

\vspace{0.2cm}

{\it Remark.} A $k$-Gorenstein ring (resp. an $\infty$-Gorenstein
ring) is clearly right quasi $k$-Gorenstein (resp. right quasi
$\infty$-Gorenstein). However, we gave examples in [H1] to explain
that there are right quasi $k$-Gorenstein rings (resp. right quasi
$\infty$-Gorenstein rings) which are not $k$-Gorenstein (resp.
$\infty$-Gorenstein); and contrary to the notion of $k$-Gorenstein
(resp. $\infty$-Gorenstein), the notion of quasi $k$-Gorenstein
(resp. quasi $\infty$-Gorenstein) is not left-right symmetric.

\vspace{0.2cm}

The following corollary develops [H1, Theorem 1].

\vspace{0.2cm}

{\bf Corollary 3.5} {\it Let} $\Lambda$ {\it be a right quasi}
$k$-{\it Gorenstein ring. Then we have}

(1) {\it Each module in} mod $\Lambda$ {\it has a}
$\mathbb{W}^{t}$-{\it approximation presentation for any} $1\leq t
\leq k+1$.

(2) {\it Each module in} mod $\Lambda ^{op}$ {\it has a}
$(\mathbb{W}^{t})^{op}$-{\it approximation presentation for any}
$1\leq t \leq k$.

\vspace{0.2cm}

{\it Proof.} (1) and (2) follow from Theorems 3.3 and 3.1,
respectively. \hfill{$\blacksquare$}

\vspace{0.2cm}

It is interesting to know when $\mathcal{P}^{k}(\Lambda)$ is
homologically finite in mod $\Lambda$. For an artin algebra
$\Lambda$ the following results are known:

(1) $\mathcal{P}^{0}(\Lambda)$ is functorially finite.

(2) $\mathcal{P}^{1}(\Lambda)$ is covariantly finite (see [AR1]);
$\mathcal{P}^{1}(\Lambda)$ is contravariantly finite if the
projective dimension of the injective envelope of $\Lambda$ as a
left $\Lambda$-module is at most one (see [IST]) (Dually, we have
that $\mathcal{P}^{1}(\Lambda ^{op})$ is contravariantly finite in
mod $\Lambda ^{op}$ if $\Lambda$ is a right quasi 1-Gorenstein
algebra).

(3) $\mathcal{P}^{k}(\Lambda)$ is functorially finite if $\Lambda$
is of finite representation type, where $k\in \mathbb{N}\bigcup
\{\infty\}$ ($\mathbb{N}$ denotes positive integers) (see [AR1]).

(4) $\mathcal{P}^{k}(\Lambda)$ is contravariantly finite if
$\Lambda$ is stably equivalent to a hereditary algebra, where $k\in
\mathbb{N}\bigcup \{\infty\}$ (see [AR1, D]).

In [AR1] Auslander and Reiten showed that if
$\mathcal{P}^{\infty}(\Lambda)$ is contravariantly finite in mod
$\Lambda$ then {\bf FDC} holds true for $\Lambda$. However,
$\mathcal{P}^{k}(\Lambda)$ and $\mathcal{P}^{\infty}(\Lambda)$ need
not to be contravariantly finite in mod $\Lambda$ when {\bf FDC}
holds true (see [IST]).

In the following we will study the homological finiteness of
$\mathcal{P}^{k}(\Lambda ^{op})$ and $\mathcal{P}^{k}(\Lambda)$ over
a right quasi $k$-Gorenstein algebra (or ring) $\Lambda$.

For a non-negative integer $n$, we use $\Omega ^{n}$(mod $\Lambda$)
to denote the full subcategory of mod $\Lambda$ consisting of
$n$-syzygy modules, and use $\mathcal{I}^{n}(\Lambda)$ to denote the
full subcategory of mod $\Lambda$ consisting of the modules with
injective dimension at most $n$.

\vspace{0.2cm}

The following result generalizes and develops [IST, Theorem 2.1] and
[AR2, Proposition 5.8].

\vspace{0.2cm}

{\bf Theorem 3.6} {\it Let} $\Lambda$ {\it be an artin right quasi}
$k$-{\it Gorenstein algebra. Then} $\mathcal{P}^{k}(\Lambda ^{op})$
{\it is functorially finite in} mod $\Lambda ^{op}$ {\it and has
almost split sequences}.

\vspace{0.2cm}

{\it Proof.} $\Omega ^{k}$(mod $\Lambda$) is functorially finite in
mod $\Lambda$ by [ASo, Section 3] and closed under extensions by
[AR3, Theorem 4.7]. So, by the dual version of [AR1, Remark after
Proposition 1.8], $\mathcal{I}^{k}(\Lambda)=\{ C\in {\rm mod}\
\Lambda |$Ext$_{\Lambda}^{1}(\Omega ^{k}({\rm mod}\ \Lambda),
C)=0\}$ is covariantly finite in mod $\Lambda$. Hence
$\mathcal{P}^{k}(\Lambda ^{op})$ is contravariantly finite in mod
$\Lambda ^{op}$. Then $\mathcal{P}^{k}(\Lambda ^{op})$ is also
covariantly finite in mod $\Lambda ^{op}$ by [KS, Corollary 2.6]. By
[ASm2, Theorem 2.4], $\mathcal{P}^{k}(\Lambda ^{op})$ has almost
split sequences. \hfill{$\blacksquare$}

\vspace{0.2cm}

It is still open when $\mathcal{P}^{k}(\Lambda)$ is covariantly (or
contravariantly) finite in mod $\Lambda$ for an (artin) right quasi
$k$-Gorenstein algebra $\Lambda$ (see [AR1]). However, as an
application of Theorem 3.3 we have the following result.

\vspace{0.2cm}

{\bf Proposition 3.7} {\it Let} $M\in$mod $\Lambda$ {\it with}
s.gradeExt$^{t+1}_{\Lambda}(M, \Lambda)\geq t$ {\it for any} $1\leq
t \leq k-1$. {\it If} $\Lambda$ {\it has property} (W$^{k}$), {\it
then} $M$ {\it has a left} $\mathcal{P}^{k}(\Lambda)$-{\it
approximation}.

\vspace{0.2cm}

{\it Proof.} Following Theorem 3.3, we assume that $0 \to K_{k}(M)
\to E_{k}(M) \to M \to 0$ is a ${\rm \mathbb{W}}^{k}$-approximation
presentation of $M$ with l.pd$_{\Lambda}K_{k}(M)\leq k-1$ and
$E_{k}(M)\in \mathbb{W}^{k}$. Because $\Lambda$ has property
(W$^{k}$), $E_{k}(M)$ is torsionless. From the proof of Theorem 3.3
we know that there is an exact sequence $0 \to E_{k}(M) \to P \to E
\to 0$ such that $E_{k}(M) \to P$ is a left
$\mathcal{P}^{0}(\Lambda)$-approximation and $E\in
\mathbb{W}^{k+1}$. Consider the following push-out diagram:

$$\xymatrix{& & 0 \ar[d] & 0 \ar[d] &\\
0 \ar[r] & K_{k}(M) \ar@{=}[d] \ar[r] & E_{k}(M) \ar[d]
\ar[r] & M \ar[d] \ar[r] & 0\\
0 \ar[r] & K_{k}(M) \ar[r] & P \ar[r] \ar[d] & F \ar[r] \ar[d] &
0\\
& & E \ar[d] \ar@{=}[r] & E \ar[d] &\\
& & 0 & 0 &}
$$
From the middle row we know that l.pd$_{\Lambda}F\leq k$. Since
$E\in \mathbb{W}^{k+1}$, it is easy to see that the exact sequence
$0 \to M \to F \to E \to 0$ is a left
$\mathcal{P}^{k}(\Lambda)$-approximation. \hfill{$\blacksquare$}

\vspace{0.2cm}

{\bf Corollary 3.8} {\it Let} $\Lambda$ {\it be a right quasi}
$k$-{\it Gorenstein ring}.

(1) {\it If} $\Lambda$ {\it has property} (W$^{k}$), {\it then}
$\mathcal{P}^{k}(\Lambda)$ {\it is covariantly finite in} mod
$\Lambda$.

(2) {\it If} $\Lambda$ {\it has property} (W$^{k+1}$), {\it then}
$\mathcal{P}^{k+1}(\Lambda)$ {\it is covariantly finite in} mod
$\Lambda$.

\vspace{0.2cm}

{\it Proof.} Notice for a right quasi $k$-Gorenstein ring $\Lambda$
we have s.gradeExt$^{t+1}_{\Lambda}(M, \Lambda)\geq t$ for any $M
\in$mod $\Lambda$ and $1\leq t \leq k$, so (1) and (2) follow from
Proposition 3.7. \hfill{$\blacksquare$}

\vspace{0.2cm}

In the rest of this section, $\Lambda$ is an artin algebra. In
addition, we give some properties for right quasi $k$-Gorenstein
algebras.

\vspace{0.2cm}

{\bf Proposition 3.9} {\it Let} $\Lambda$ {\it be a right quasi}
$k$-{\it Gorenstein ring. Then} $\mathbb{W}^{k}\subseteq \{ M\in
{\rm mod}\ \Lambda |$ \linebreak ${\rm Ext}_{\Lambda}^{1}(M,
\mathcal{P}^{k-1}(\Lambda))=0\}$ {\it and} $(\mathbb{W}^{k})^{op}=\{
N\in {\rm mod}\ \Lambda ^{op}|{\rm Ext}_{\Lambda}^{1}(N,
\mathcal{P}^{k-1}(\Lambda ^{op}))=0\}$.

\vspace{0.2cm}

{\it Proof.} It is trivial that $\mathbb{W}^{k}\subseteq \{ M\in
{\rm mod}\ \Lambda |{\rm Ext}_{\Lambda}^{1}(M,
\mathcal{P}^{k-1}(\Lambda))=0\}$ and $(\mathbb{W}^{k})^{op}\subseteq
\{ N\in {\rm mod}\ \Lambda ^{op}|{\rm Ext}_{\Lambda}^{1}(N,
\mathcal{P}^{k-1}(\Lambda ^{op}))=0\}$. Now, to show the inclusion
$(\mathbb{W}^{k})^{op}\supseteq \{ N\in {\rm mod}\ \Lambda
^{op}|{\rm Ext}_{\Lambda}^{1}(N,$ \linebreak
$\mathcal{P}^{k-1}(\Lambda ^{op}))=0\}$, consider the minimal
injective resolution of $\Lambda$ as a right $\Lambda$-module: $0\to
\Lambda _{\Lambda}\to I_{0} \to I_{1} \to \cdots \to I_{i} \to
\cdots$. Because $\Lambda$ is an artin right quasi $k$-Gorenstein
ring, l.pd$_{\Lambda}I_{i}\leq i+1$ for any $0\leq i \leq k-1$. Then
inductively we have l.pd$_{\Lambda}{\rm Im}(I_{i}\to I_{i+1})\leq
i+1$ for any $0\leq i \leq k-1$. Now it is not difficult to verify
that if a module $N\in$ mod $\Lambda ^{op}$ satisfies ${\rm
Ext}_{\Lambda}^{1}(N, \mathcal{P}^{k-1}(\Lambda ^{op}))=0$ then
$N\in(\mathbb{W}^{k})^{op}$. \hfill{$\blacksquare$}

\vspace{0.2cm}

Since the notion of $k$-Gorenstein rings is left-right aymmetric, by
Proposition 3.9 we have the following

\vspace{0.2cm}

{\bf Corollary 3.10} {\it Let} $\Lambda$ {\it be a} $k$-{\it
Gorenstein ring. Then} $\mathbb{W}^{k}= \{ M\in {\rm mod}\ \Lambda
|{\rm Ext}_{\Lambda}^{1}(M, \mathcal{P}^{k-1}(\Lambda))$ \linebreak
$=0\}$ {\it and} $(\mathbb{W}^{k})^{op}=\{ N\in {\rm mod}\ \Lambda
^{op}|{\rm Ext}_{\Lambda}^{1}(N, \mathcal{P}^{k-1}(\Lambda
^{op}))=0\}$.

\vspace{0.2cm}

{\bf Theorem 3.11} {\it Let} $\Lambda$ {\it be a right quasi}
$k$-{\it Gorenstein ring. Then for any} $N\in$ mod $\Lambda ^{op}$
{\it there is an exact sequence}:
$$N \to K^{(k)}(N)\to E^{(k)}(N)\to 0 \eqno{\textcircled{\small{5}}}$$
{\it where} $N \to K^{(k)}(N)$ {\it is a minimal left}
$\mathcal{P}^{k-1}(\Lambda ^{op})$-{\it approximation of} $N$ {\it
and} $E^{(k)}(N)\in (\mathbb{W}^{k})^{op}$.

\vspace{0.2cm}

{\it Proof.} By Theorem 3.6, $\mathcal{P}^{k-1}(\Lambda ^{op})$ is
covariantly finite in mod $\Lambda ^{op}$ and any $N\in$ mod
$\Lambda ^{op}$ has a minimal left $\mathcal{P}^{k-1}(\Lambda
^{op})$-approximation $N \to K^{(k)}(N)$. Put $E^{(k)}(N)=$Coker$(N
\to K^{(k)}(N))$. Then by Wakamatsu's Lemma (see [AR1, Lemma 1.3]),
we have ${\rm Ext}_{\Lambda}^{1}(E^{(k)}(N),
\mathcal{P}^{k-1}(\Lambda ^{op}))=0$. So $E^{(k)}(N) \in
(\mathbb{W}^{k})^{op}$ by Proposition 3.9 and we get the exact
sequence \textcircled{\small{5}}. \hfill{$\blacksquare$}

\vspace{0.2cm}

A natural question is: when is the minimal left
$\mathcal{P}^{k-1}(\Lambda ^{op})$-approximation of $N$: $N \to
K^{(k)}(N)$ in the exact sequence \textcircled{\small{5}}
monomorphic? If $N$ is torsionless in mod $\Lambda ^{op}$, then the
answer to this question is clearly affirmative.

\vspace{0.5cm}

\centerline{\large \bf 4. Property (W$^{k}$) and homological
conjectures}

\vspace{0.2cm}

In [AR2] Auslander and Reiten conjecture that an artin
$\infty$-Gorenstein algebra is Gorenstein ({\bf ARC}). The famous
Nakayama Conjecture ({\bf NC}) is a special case of {\bf FDC} (resp.
{\bf ARC}). In this section we give a connection between {\bf FDC},
{\bf ARC} and property (W$^{k})$ as follows.

\vspace{0.2cm}

{\bf Theorem 4.1} {\it Let} $k$ {\it be a non-negative integer. The
following statements are equivalent for an artin right quasi
$\infty$-Gorenstein algebra} $\Lambda$.

(1) $\Lambda$ {\it is Gorenstein with self-injective dimension at
most} $k$.

(2) l.id$_{\Lambda}\Lambda \leq k$.

(2)$^{op}$ r.id$_{\Lambda}\Lambda \leq k$.

(3) $\Lambda$ {\it has property} (W$^{k})$.

(3)$^{op}$ $\Lambda$ {\it has property} (W$^{k})^{op}$.

{\it In particular, we have} ${\rm fin.dim}\Lambda \leq {\rm
l.id}_{\Lambda}\Lambda \leq {\rm fin.dim}\Lambda+1$.

\vspace{0.2cm}

In [AR2, Corollary 5.5], Auslander and Reiten showed that the
conditions (1) and (2) in Theorem 4.1 are equivalent for an artin
$\infty$-Gorenstein algebra (also see [Be, Corollary 6.21]). Notice
that {\bf NC} is a special case of {\bf ARC}, so we know from
Theorem 4.1 that in order to verify {\bf NC} (for any artin algebra)
it suffices to verify {\bf FDC} for (right quasi)
$\infty$-Gorenstein algebras.

To prove this theorem we need some lemmas.

\vspace{0.2cm}

{\bf Lemma 4.2} {\it Let} $\Lambda$ {\it be a right quasi} $k$-{\it
Gorenstein ring (especially, a right quasi $\infty$-Gorenstein ring)
with property} (W$^{k})^{op}$. {\it Then} l.id$_{\Lambda}\Lambda
\leq k$.

\vspace{0.2cm}

{\it Proof.} Let $\Lambda$ be a right quasi $k$-Gorenstein ring and
$M$ any module in mod $\Lambda$. By Corollary 3.5(1) there is a
${\rm \mathbb{W}}^{k}$-approximation presentation of $M$:
$$0\to K_{k}(M) \to E_{k}(M) \to M \to 0$$
with l.pd$_{\Lambda}K_{k}(M) \leq k-1$ and $E_{k} \in {\rm
\mathbb{W}}^{k}$. Since $\Lambda$ has property (W$^{k})^{op}$, it
follows from the dual version of Lemma 2.1 that $E_{k}$ is a
W$^{\infty}$-module and we then have Ext$_{\Lambda}^{i+1}(M,
\Lambda)\cong$Ext$_{\Lambda}^{i}(K_{k}(M), \Lambda)=0$ for any
$i\geq k$, which implies that l.id$_{\Lambda}\Lambda \leq k$.
\hfill{$\blacksquare$}

\vspace{0.2cm}

By using Corollary 3.5(2) and an argument similar to that in proving
Lemma 4.2, we have the following result, which says that the
converse of Corollary 2.3(1) holds true for a right quasi
$k$-Gorenstein ring $\Lambda$.

\vspace{0.2cm}

{\bf Lemma 4.3} {\it Let} $\Lambda$ {\it be a right quasi} $k$-{\it
Gorenstein ring (especially, a right quasi $\infty$-Gorenstein ring)
with property} (W$^{k})$. {\it Then} r.id$_{\Lambda}\Lambda \leq k$.

\vspace{0.2cm}

{\bf Lemma 4.4} ([H3, Corollary 4]) {\it If} $\Lambda$ {\it is an
artin right quasi $\infty$-Gorenstein algebra, then}
l.id$_{\Lambda}\Lambda =$ \linebreak r.id$_{\Lambda}\Lambda$.

\vspace{0.2cm}

The following result is well known.

\vspace{0.2cm}

{\bf Lemma 4.5} l.id$_{\Lambda}\Lambda \geq$fin.dim$\Lambda$.

\vspace{0.2cm}

{\bf Lemma 4.6} {\it Let} $\Lambda$ {\it be a right quasi
$\infty$-Gorenstein ring. Then} l.id$_{\Lambda}\Lambda \leq
$fin.dim$\Lambda +1$.

\vspace{0.2cm}

{\it Proof.} Without loss of generality, we assume that
fin.dim$\Lambda =k<\infty$ and $M$ is any module in mod $\Lambda$.
By Corollary 3.5(1), for any $i\geq 1$ there is a ${\rm
\mathbb{W}}^{k+i+1}$-approximation presentation of $M$:
$$0\to K_{k+i+1}(M) \to E_{k+i+1}(M) \to M \to 0$$
with l.pd$_{\Lambda}K_{k+i+1}(M) \leq k+i$ and $E_{k+i+1} \in {\rm
\mathbb{W}}^{k+i+1}$. Then l.pd$_{\Lambda}K_{k+i+1}(M) \leq k$ since
fin.dim$\Lambda =k$. So Ext$_{\Lambda}^{k+i+1}(M,
\Lambda)\cong$Ext$_{\Lambda}^{k+i}(K_{k+i+1}(M), \Lambda)=0$ and
l.id$_{\Lambda}\Lambda \leq k+1$. We are done.
\hfill{$\blacksquare$}

\vspace{0.2cm}

{\it Proof of Theorem 4.1.} By Lemma 4.4 we get the equivalence of
(1), (2) and (2)$^{op}$. That (2) (resp. (2)$^{op}$) implies
(3)$^{op}$ (resp. (3)) follows from the dual version of [HT, Theorem
2.2] (resp. Corollary 2.3), and (3)$^{op}$ (resp. (3)) implies (2)
(resp. (2)$^{op}$) follows from Lemma 4.2 (resp. Lemma 4.3). The
last assertion follows from Lemmas 4.5 and 4.6.
\hfill{$\blacksquare$}

\vspace{0.2cm}

{\bf Corollary 4.7} {\it Let} $\Lambda$ {\it be an artin right quasi
$\infty$-Gorenstein algebra. If} $\Lambda _{\Lambda}$ {\it has a
ultimately closed injective resolution, then} $\Lambda$ {\it is
Gorenstein}.

\vspace{0.2cm}

{\it Proof.} By Theorem 4.1 and Corollary 2.3.
\hfill{$\blacksquare$}

\vspace{0.5cm}

\centerline{\large \bf 5. Finite homological dimensions}

\vspace{0.2cm}

Let $M\in$ mod $\Lambda$ and $n$ a non-negative integer. Recall that
$M$ is said to have {\it Gorenstein dimension zero} if it satisfies
the conditions: (1) $M$ is reflexive; (2) $M\in
{\mathbb{W}^{\infty}}$ and $M^{*}\in ({\mathbb{W}^{\infty}})^{op}$;
and $M$ is said to have {\it finite Gorenstein dimension} $n$,
written as G-dim$_{\Lambda}M=n$, if $n$ is the least non-negative
integer such that there is an exact sequence $0 \to X_{n} \to \cdots
\to X_{1} \to X_{0} \to M \to 0$ with all $X_{i}$ having Gorenstein
dimension zero (see [AB]). Gorenstein dimension of modules is an
important invariant in homological algebra (see [C]). It is clear
that in general
$\mathbb{W}^{\infty}$-dim$_{\Lambda}M\leq$G-dim$_{\Lambda}M$. In
[H2] we showed that for an artin algebra $\Lambda$,
l.id$_{\Lambda}\Lambda <\infty$ if and only if
$\mathbb{W}^{\infty}$-dim$_{\Lambda}M<\infty$ for each $M\in$ mod
$\Lambda$, if and only if G-dim$_{\Lambda}M<\infty$ for each $M\in$
mod $\Lambda$. Here we give the following

\vspace{0.2cm}

{\bf Proposition 5.1} {\it If} $\Lambda$ {\it has property}
(W$^{\infty}$), {\it then for any} $M\in$ mod $\Lambda$ {\it we
have} $\mathbb{W}^{\infty}$-dim$_{\Lambda}M=$G-dim$_{\Lambda}M$.

\vspace{0.2cm}

{\it Proof.} For any $M \in$mod $\Lambda$, it is enough to show that
G-dim$_{\Lambda}M \leq {\mathbb{W}^{\infty}}$-dim$_{\Lambda}M$ in
case $\mathbb{W}^{\infty}$-dim$_{\Lambda}M=k<\infty$, because
$\mathbb{W}^{\infty}$-dim$_{\Lambda}M \leq $G-dim$_{\Lambda}M$ by
definition. We prove it by induction on $k$.

If $k=0$, then $M$ is in $\mathbb{W}^{\infty}$ and by assumption $M$
is torsionless. Let $P \buildrel {f} \over \to M^{*} \to 0$ be exact
in mod $\Lambda ^{op}$ with $P$ projective. Then by Lemma 3.2, the
composition: $M \buildrel {\sigma _{M}} \over \hookrightarrow M^{**}
\buildrel {f^{*}} \over \hookrightarrow P^{*}$ is a left
$\mathcal{P}^{0}(\Lambda)$-approximation of $M$. Put $g=f^{*}\sigma
_{M}$ and $N=$Coker$g$. Then Ext$_{\Lambda}^{1}(N, \Lambda)=0$. Thus
from the exactness of $0\to M \to P^{*} \to N \to 0$ we have the
following commutative diagram with exact rows:
$$\xymatrix{0\ar[r] & M \ar[r] \ar[d]^{\sigma _{M}}
& P^{*} \ar[r] \ar[d]^{\sigma _{P^{*}}} & N \ar[r] \ar[d]^{\sigma
_{N}} & 0 \\
0\ar[r] & M^{**} \ar[r] & P^{***} \ar[r] & N^{**} }$$ where $\sigma
_{M}$ is a monomorphism and $\sigma _{P^{*}}$ is an isomorphism.
Since $M\in {\mathbb{W}^{\infty}}$, it is easy to see that $N\in
{\mathbb{W}^{\infty}}$ and by assumption $N$ is also torsionless.
Thus $\sigma _{N}$ is a monomorphism and therefore $\sigma _{M}$ is
an isomorphism and $M$ is reflexive. Similarly we know that $N$ is
also reflexive. So, also from the exact sequence $0\to M \to P^{*}
\to N \to 0$ we have an exact sequence $0\to M^{**} \to P^{***} \to
N^{**} \to 0$ in mod $\Lambda$. Notice that $0\to N^{*} \to P^{**}
\to M^{*} \to 0$ is exact in mod $\Lambda ^{op}$, and
Ext$_{\Lambda}^{1}(M^{*}, \Lambda)=0$. Similar to the above argument
we have Ext$_{\Lambda}^{1}(N^{*}, \Lambda)=0$ and so
Ext$_{\Lambda}^{2}(M^{*}, \Lambda)=0$. Continuing this process, we
finally get that $M^{*}\in (\mathbb{W}^{\infty})^{op}$ and
G-dim$_{\Lambda}M=0$.

Now suppose $k\geq 1$ and $0 \to L \to P_{k-1} \to \cdots \to P_{0}
\to M \to 0$ is exact in mod $\Lambda$ with all $P_{i}$ projective.
Then $\mathbb{W}^{\infty}$-dim$_{\Lambda}L=0$ and by the argument
above we have G-dim$_{\Lambda}L=0$. It follows from [AB, Theorem
3.13] that G-dim$_{\Lambda}M\leq k$. This completes the proof.
\hfill{$\blacksquare$}

\vspace{0.2cm}

Assume that $\Lambda$ has property (W$^{\infty}$). If $M$ is a
W$^{\infty}$-module in mod $\Lambda$, then
$\mathbb{W}^{\infty}$-dim$_{\Lambda}M=0$ and by Proposition 5.1 we
have G-dim$_{\Lambda}M=0$. Thus $M$ is reflexive and consequently we
get the following interesting result (compare it with Lemma 2.1).

\vspace{0.2cm}

{\bf Corollary 5.2} {\it The following statements are equivalent}.

(1) $\Lambda$ {\it has property} (W$^{\infty}$), {\it that is}, {\it
each} ${\rm W}^{\infty}$-{\it module in} mod $\Lambda$ {\it is
torsionless}.

(2) {\it Each} ${\rm W}^{\infty}$-{\it module in} mod $\Lambda$
{\it is reflexive}.

\vspace{0.2cm}

Let $\mathcal{X}$ be a contravariantly finite subcategory of mod
$\Lambda$. Then for each module $M$ in mod $\Lambda$ we have a
complex:
$$\cdots \buildrel {f_{n+1}} \over \longrightarrow
X_{n}\buildrel {f_{n}} \over \to \cdots \buildrel {f_{2}} \over \to
X_{1}\buildrel {f_{1}} \over \to X_{0} \buildrel {f_{0}} \over \to
M\to 0 \eqno{\textcircled{\small{6}}}$$ with $X_{0}\to M$ a right
$\mathcal{X}$-approximation of $M$ and each $X_{i}\to {\rm
Ker}f_{i-1}$ a right $\mathcal{X}$-approximation of ${\rm
Ker}f_{i-1}$ for any $i\geq 1$. Since $\mathcal{X}$ is
contravariantly finite, we have an exact sequence:
$$\cdots \to {\rm Hom}_{\Lambda}(\mathcal{X}, X_{n}) \to \cdots \to
{\rm Hom}_{\Lambda}(\mathcal{X}, X_{1}) \to {\rm
Hom}_{\Lambda}(\mathcal{X}, X_{0}) \to {\rm
Hom}_{\Lambda}(\mathcal{X}, M) \to 0$$ We define pd$_{(\mathcal{X},
{\rm mod\ }\Lambda)}M$ to be inf$\{ n | 0 \to {\rm
Hom}_{\Lambda}(\mathcal{X}, X_{n}) \to \cdots \to {\rm
Hom}_{\Lambda}(\mathcal{X}, X_{1}) \to {\rm
Hom}_{\Lambda}(\mathcal{X}, X_{0}) \to {\rm
Hom}_{\Lambda}(\mathcal{X}, M) \to 0$ is exact$\}$. If no such $n$
exists set pd$_{(\mathcal{X}, {\rm mod\ }\Lambda)}M=\infty$. We also
define gl.dim$(\mathcal{X}, {\rm mod\ }\Lambda)$ to be sup$\{$
pd$_{(\mathcal{X}, {\rm mod\ }\Lambda)}M|M\in$mod $\Lambda \}$.

On the other hand, if there is some $n$ such that the complex
\textcircled{\small{6}} stops after $n$ steps, that is, we have a
complex $0\to X_{n}\buildrel {f_{n}} \over \to \cdots \buildrel
{f_{2}} \over \to X_{1}\buildrel {f_{1}} \over \to X_{0} \buildrel
{f_{0}} \over \to M\to 0$ with $X_{0}\to M$ a right
$\mathcal{X}$-approximation of $M$ and each $X_{i}\to {\rm
Ker}f_{i-1}$ a right $\mathcal{X}$-approximation of ${\rm
Ker}f_{i-1}$ for any $1\leq i \leq n$, then pd$_{\mathcal{X}}M$ is
defined to be the smallest non-negative $n$ such that the complex
\textcircled{\small{6}} stops after $n$ steps. If no such $n$ exists
set pd$_{\mathcal{X}}M=\infty$. In addition,
gl.dim$_{\mathcal{X}}\Lambda$ is defined to be
sup$\{$pd$_{\mathcal{X}}M|M\in$mod $\Lambda \}$.

\vspace{0.2cm}

{\bf Proposition 5.3} {\it Let} $\mathcal{X}$ {\it be a
contravariantly finite subcategory of} mod $\Lambda$ {\it
containing} $\mathcal{P}^{0}(\Lambda)$. {\it Then}
pd$_{(\mathcal{X}, {\rm mod\ }\Lambda)}M=$pd$_{\mathcal{X}}M$ {\it
for any} $M \in$ mod $\Lambda$.

\vspace{0.2cm}

{\it Proof.} Clearly we have pd$_{(\mathcal{X}, {\rm mod\
}\Lambda)}M\leq$pd$_{\mathcal{X}}M$. Now suppose pd$_{(\mathcal{X},
{\rm mod\ }\Lambda)}M=n<\infty$. Because $\mathcal{P}^{0}(\Lambda)$
is contained in $\mathcal{X}$, the complex \textcircled{\small{6}}
is exact. By definition of pd$_{(\mathcal{X}, {\rm mod\
}\Lambda)}M$,
$$0 \to {\rm Hom}_{\Lambda}(\mathcal{X}, X_{n}) \to
{\rm Hom}_{\Lambda}(\mathcal{X}, X_{n-1}) \to \cdots \to {\rm
Hom}_{\Lambda}(\mathcal{X}, X_{1}) \to $$ $${\rm
Hom}_{\Lambda}(\mathcal{X}, X_{0}) \to {\rm
Hom}_{\Lambda}(\mathcal{X}, M) \to 0$$ is also exact. On the other
hand, if we put $K_{i}=$Ker$f_{i-1}$ for any $i\geq 1$, we also have
an exact sequence:
$$0 \to {\rm Hom}_{\Lambda}(\mathcal{X}, K_{n}) \to
{\rm Hom}_{\Lambda}(\mathcal{X}, X_{n-1}) \to \cdots \to {\rm
Hom}_{\Lambda}(\mathcal{X}, X_{1}) \to $$ $${\rm
Hom}_{\Lambda}(\mathcal{X}, X_{0}) \to {\rm
Hom}_{\Lambda}(\mathcal{X}, M) \to 0.$$ So ${\rm
Hom}_{\Lambda}(\mathcal{X}, X_{n})\cong {\rm
Hom}_{\Lambda}(\mathcal{X}, K_{n})$ and hence ${\rm
Hom}_{\Lambda}(\mathcal{X}, K_{n+1})=0$. Now consider the exact
sequence $0\to K_{n+2} \buildrel {\alpha} \over \to X_{n+1} \to
K_{n+1} \to 0$. We then have that ${\rm Hom}_{\Lambda}(\mathcal{X},
K_{n+2}) \buildrel {{\rm Hom}_{\Lambda}(\mathcal{X}, \alpha}) \over
\longrightarrow{\rm Hom}_{\Lambda}(\mathcal{X}, X_{n+1})$ is an
isomorphism, which yields in particular that ${\rm
Hom}_{\Lambda}(X_{n+1}, K_{n+2}) \buildrel {{\rm
Hom}_{\Lambda}(X_{n+1}, \alpha}) \over \longrightarrow{\rm
Hom}_{\Lambda}(X_{n+1}, X_{n+1})$ is an isomorphism and $\alpha$ is
an epimorphism and hence an isomorphism. So $K_{n+1}=0$ and
pd$_{\mathcal{X}}M\leq n$. This finishes the proof.
\hfill{$\blacksquare$}

\vspace{0.2cm}

{\bf Corollary 5.4} {\it Under the assumptions of Proposition 5.3,
we have} gl.dim$(\mathcal{X}, {\rm mod\
}\Lambda)=$gl.dim$_{\mathcal{X}}\Lambda$.

\vspace{0.2cm}

The following corollary is an immediate consequence of Corollary
3.5.

\vspace{0.2cm}

{\bf Proposition 5.5} {\it Let} $\Lambda$ {\it be a right quasi}
$k$-{\it Gorenstein ring. Then for any} $M\in$mod $\Lambda$ {\it
and} $1\leq t \leq k+1$, {\it we have} pd$_{(\mathbb{W}^{t}, {\rm
mod\ }\Lambda)}M$=pd$_{\mathbb{W}^{t}}M\leq t$ {\it and so}
gl.dim$(\mathbb{W}^{t}, {\rm mod\
}\Lambda)$=gl.dim$_{\mathbb{W}^{t}}\Lambda \leq t$.

\vspace{0.2cm}

{\it Proof.} By Corollary 3.5(1), for any $1\leq t \leq k+1$,
$\mathbb{W}^{t}$ is a contravariantly finite subcategory of mod
$\Lambda$ containing $\mathcal{P}^{0}(\Lambda)$. Then by Proposition
5.3 and Corollary 5.4, we have that pd$_{(\mathbb{W}^{t}, {\rm mod\
}\Lambda)}M$=pd$_{\mathbb{W}^{t}}M$ for any $M\in$mod $\Lambda$ and
gl.dim$(\mathbb{W}^{t}, {\rm mod\
}\Lambda)$=gl.dim$_{\mathbb{W}^{t}}\Lambda$.

By Corollary 3.5(1), for any $M\in$mod $\Lambda$ and $1\leq t \leq
k+1$ there is a $\mathbb{W}^{t}$-approximation presentation of $M$:
$$0 \to K_{t}(M) \to E_{t}(M) \to M \to 0$$ with
l.pd$_{\Lambda}K_{t}(M)\leq t-1$ and $E_{t}(M)\in \mathbb{W}^{t}$.

For a W$^{t}$-module $A$, it is easy to see that $A$ is projective
if l.pd$_{\Lambda}A\leq t-1$. Now, considering the
$\mathbb{W}^{t}$-approximation presentation of $K_{t}(M)$ we get
easily our conclusion. \hfill{$\blacksquare$}

\vspace{0.2cm}

For a subcategory $\mathcal{X}$ of mod $\Lambda$ we use
$\hat{\mathcal{X}}$ to denote the subcategory of mod $\Lambda$
consisting of the modules $C$ for which there is an exact sequence
$0 \to X_{n} \to \cdots \to X_{1} \to X_{0} \to C \to 0$ with all
$X_{i}$ in $\mathcal{X}$ (see [AR1]).

\vspace{0.2cm}

{\bf Theorem 5.6} {\it Let} $\Lambda$ {\it be an artin right quasi
$\infty$-Gorenstein algebra}. {\it Then the following statements are
equivalent}.

(1) $\Lambda$ {\it is Gorenstein}.

(2) $\mathbb{W}^{\infty}$ {\it is contravariantly finite and}
gl.dim$(\mathbb{W}^{\infty}, {\rm mod\ }\Lambda)$ {\it is finite}.

(3) $\mathbb{W}^{\infty}$ {\it is contravariantly finite and}
pd$_{(\mathbb{W}^{\infty}, {\rm mod\ }\Lambda)}M$ {\it is finite
for any} $M\in$mod $\Lambda$.

(4) $\mathbb{W}^{\infty}$ {\it is contravariantly finite and}
gl.dim$_{\mathbb{W}^{\infty}}\Lambda$ {\it is finite}.

(5) $\mathbb{W}^{\infty}$ {\it is contravariantly finite and}
pd$_{\mathbb{W}^{\infty}}M$ {\it is finite for any} $M\in$mod
$\Lambda$.

(6) fin.dim$\Lambda$ {\it is finite}.

\vspace{0.2cm}

{\it Proof.} The equivalence of (1) and (6) follows from Theorem
4.1, and that (2) implies (3) and (4) implies (5) are trivial.
Because $\mathcal{P}^{0}(\Lambda)$ is contained in
$\mathbb{W}^{\infty}$, from Corollary 5.4 (resp. Proposition 5.3) we
know that (2) and (4) (resp. (3) and (5)) are equivalent. Assume
that $\Lambda$ is Gorenstein. Then
$\mathbb{W}^{\infty}=\mathbb{W}^{t}$ for some positive integer $t$
and from Proposition 5.5 we know that $(1) \Rightarrow (2)$ holds.
Now suppose (5) holds. Then we have that
$\widehat{\mathbb{W}^{\infty}}=$mod $\Lambda$. Moreover,
$\mathbb{W}^{\infty}$ is clearly resolving (that is,
$\mathbb{W}^{\infty}$ is closed under extensions, kernels of
epimorphisms, and contains $\mathcal{P}^{0}(\Lambda)$). It follows
from [AR1, Theorem 5.5] that $_{\Lambda}\Lambda$ is cotilting and so
$\Lambda$ is Gorenstein. This proves (5) implies (1). We are done.
\hfill{$\blacksquare$}

\vspace{0.2cm}

In view of Theorem 5.6 it would be interesting to know when
$\mathbb{W}^{\infty}$ is contravariantly finite.

\vspace{0.5cm}

\centerline{\large \bf 6. Conjectures}

\vspace{0.2cm}

From the result obtained above, as far as $\mathbb{W}^{\infty}$ is
concerned, two problems are worth being studied: when is a module in
$\mathbb{W}^{\infty}$ torsionless (that is, when has a ring
$\Lambda$ property (W$^{\infty}$)) and what does property
(W$^{\infty}$) imply about $\Lambda$? In particular, we pose the
following

\vspace{0.2cm}

 {\bf
Conjecture I} Over an artin algebra, each module in
$\mathbb{W}^{\infty}$ is torsionless (equivalently, each module in
$\mathbb{W}^{\infty}$ is reflexive), that is, any artin algebra has
property (W$^{\infty})$.

\vspace{0.2cm}

{\it Remark.} (1) We know from Proposition 5.1 that if a left and
right noetherian ring $\Lambda$ has property (W$^{\infty}$) then for
any $M\in$ mod $\Lambda$ we have
$\mathbb{W}^{\infty}$-dim$_{\Lambda}M=$G-dim$_{\Lambda}M$. So, we
also conjecture that over an artin algebra the notion of the left
orthogonal dimension and that of the Gorenstein dimension of any
module in mod $\Lambda$ coincide. This is an immediate corollary of
Conjecture I.

(2) Assume that Conjecture I holds true for $\Lambda$. If
l.id$_{\Lambda}\Lambda \leq k$, then every W$^{k}$-module in mod
$\Lambda$ is a W$^{\infty}$-module and hence it is reflexive
(compare this result with Theorem 4.1).

(3) Assume that Conjecture I holds true for an artin right quasi
$k$-Gorenstein algebra $\Lambda$. If $\Lambda$ has property
(W$^{k})^{op}$, then l.id$_{\Lambda}\Lambda \leq k$ by Lemma 4.2 and
each W$^{k}$-module in mod $\Lambda$ is a W$^{\infty}$-module. By
assumption (Conjecture I holds true for $\Lambda$), we know that
each W$^{k}$-module in mod $\Lambda$ is torsionless. Similar to the
argument of Proposition 3.7, we then have that every module in mod
$\Lambda$ has a left $\mathcal{P}^{k}(\Lambda)$-approximation and
$\mathcal{P}^{k}(\Lambda)$ is covariantly finite in mod $\Lambda$
(compare this result with Corollary 3.8).

(4) The validity of Conjecture I would imply the validity of {\bf
GSC} (see [AHT, Theorem 3.10]).

\vspace{0.2cm}

In view of the results obtained in Section 4 and Theorem 5.6 we pose
the following

\vspace{0.1cm}

{\bf Conjecture II} (Generalized Auslander-Reiten Conjecture) An
artin algebra is Gorenstein if it is right quasi
$\infty$-Gorenstein.

 \vspace{0.2cm}

{\bf Acknowledgements} The research of the author was partially
supported by Specialized Research Fund for the Doctoral Program of
Higher Education (Grant No. 20060284002) and NSF of Jiangsu Province
of China (Grant No. BK2005207). The author thanks the referee for
the useful and detailed suggestions on this paper.

\end{document}